\documentclass[12pt,leqno]{article}
\usepackage{amsfonts,amsthm,amsmath}

\usepackage[mathscr]{eucal}
\usepackage{amsmath}
\usepackage{amsfonts}
\usepackage{amssymb}
\usepackage{footnpag}
\usepackage[dvips]{graphicx}
\usepackage{lscape} 

\theoremstyle{plain}
\newtheorem{thm}{Theorem}[section]

\theoremstyle{definition}
\newtheorem{df}{Definition}[section]
\newtheorem{rem}{Remark}[section]

\newtheorem{prob}{Problem}[section]

\newcommand{\Q}{\mathbb{Q}}

\newcommand{\R}{\mathbb{R}}

\newcommand{\QQ}{\mathbb{Q}}

\newcommand{\ZZ}{\mathbb{Z}}
\newcommand{\RR}{\mathbb{R}}

\newcommand{\NN}{\mathbb{N}}

\newcommand{\Cl}{\mathop{\mathrm{Cl}}\nolimits}

\begin{document}

\title{New invariants for integral lattices
}

\author
{
Ryota Hayasaka\thanks{
Drecom Co., Ltd.
ryota.hayasaka@drecom.co.jp}, 
Tsuyoshi Miezaki
\thanks{Faculty of Education, University of the Ryukyus, Okinawa  
903--0213, Japan 
miezaki@edu.u-ryukyu.ac.jp
}, and 
Masahiko Toki
\thanks{Oita National College of Technology, 
toki@oita-ct.ac.jp
}
}
\date{}
\maketitle
\begin{abstract}
Let $\Lambda$ be any integral lattice in Euclidean space. 
It has been shown that for every integer $n>0$, 
there is a hypersphere 
that passes through exactly $n$ points of $\Lambda$. 
Using this result, we introduce new lattice invariants 
and give some computational 
results related to two-dimensional Euclidean lattices 
of class number one. 
\end{abstract}
\vspace{-0.15in}

\quad

\section{Introduction}

We consider the following condition 
on lattices $\Lambda \subset \RR^{d}$. 
\begin{df}[\cite{M2,BM}]
If there is a hypersphere in $\R^{d}$ that 
passes through exactly $n$ points of $\Lambda$ for every integer $n>0$, 
then $\Lambda$ is called ``universally concyclic." 
\end{df}

A lattice generated by $(a, b), (c, d)\in \RR^{2}$, 
$(ad-bc\neq 0)$, is denoted by $\Lambda[(a, b), (c, d)]$. 
In \cite{M2}, Maehara introduced the term ``universally concyclic." 
Then, he and others showed the following results. 
In \cite{S} and \cite{MM}, Schinzel, Maehara, and Matsumoto proved that 
$\ZZ^{2}$, that is, $\Lambda[(1, 0), (0, 1)]$, is universally concyclic. 
Moreover, if $a,b,c,d\in\ZZ$ are such that $q:=ad-bc$ is a prime 
and $q\equiv 3 \pmod{4}$, 
then $\Lambda[(a, b), (c, d)]$ is universally concyclic. 
The equilateral triangular lattice 
$\Lambda[(1, 0), (-1/2, \sqrt{3}/2])]$
and 
the rectangular lattice 
$\Lambda[(1, 0), (0, \sqrt{3}])]$ 
are universally concyclic. 
In \cite{BM}, it was shown that 
all integral lattices in $\RR^d$ with $d\geq 2$ are universally concyclic. 

\begin{rem}
We remark that there exist some nonintegral lattices 
that are not universally concyclic. 
Maehara also proved in \cite{M2} that 
if $\tau$ is a transcendental number, then $\Lambda[(1, \tau), (0, 1)]$ 
cannot contain four concyclic points, and hence, it is not universally concyclic. 
The rectangular lattice $\Lambda[(\alpha, 0), (0, \beta)]$ does not contain 
five concyclic points if and only if $(\alpha /\beta )^{2}$ is an irrational number. 
Hence, some additional integrality 
conditions are necessary to ensure this property. 
\end{rem}

Let $K=\Q(\sqrt{-d})$ be an imaginary quadratic field, 
and let $\mathcal{O}_{K}$ be its ring of algebraic integers. 
Let $\Cl_{K}$ be the ideal classes of $K$. In this paper, 
we only consider the cases $\vert \Cl_{K}\vert =1$, namely, 
$d$ is in the following set: $\{1,2,3, 7, 11, 19, 43, 67,
163\}$. 

We denote by $d_{K}$ the discriminant of $K$:
\begin{eqnarray*}
d_{K}=\left\{
\begin{array}{lll}
-4d\ &{\rm if }\ -d\equiv 2,\ 3&\pmod {4} \\
-d\ &{\rm if }\ -d\equiv 1 &\pmod {4}. 
\end{array}
\right.
\end{eqnarray*}
\begin{thm}[cf.~{\cite[p.\ 87]{Zagier}}]\label{thm:lattice}
Let $d$ be a positive square-free integer, and let $K=\Q(\sqrt{-d})$. 
Then 
\begin{eqnarray*}
\mathcal{O}_{K}=
\left\{
\begin{array}{lll}
\ZZ+\ZZ\,\sqrt{-d}\quad &if\ -d\equiv 2,\ 3 &\pmod {4} \\
\ZZ+\ZZ\,\displaystyle\frac{-1+\sqrt{-d}}{2}\quad &if\ -d\equiv 1 &\pmod {4}. 
\end{array}
\right.
\end{eqnarray*}
\end{thm}
\vspace{20pt}
Therefore, we consider $\mathcal{O}_{K}$ to be a lattice in $\R^{2}$ 
with the basis 
\begin{eqnarray*}
\left\{
\begin{array}{lll}
\displaystyle(1, 0), (0,\sqrt{d} )\quad &\mbox{if}\ -d\equiv 2,\ 3 &\pmod {4} \\
\displaystyle(1, 0), \Big(-\frac{1}{2},\frac{\sqrt{d}}{2}\Big)\quad &
\mbox{if}\ -d\equiv 1 &\pmod {4}, 
\end{array}
\right.
\end{eqnarray*}
denoted by $[1,\sqrt{-d}]$, $[1,(-1+\sqrt{-d})/2]$, respectively. 
Note that $[1,\sqrt{-1}]$ is the $\ZZ^2$ lattice. 

The main purpose of this paper is 
to introduce the new lattice invariants (Definition \ref{def:main})
and to give some computational 
results related to two-dimensional Euclidean lattices of class number one
(Theorem \ref{thm:main}). 

We introduce the following new lattice invariants $\mbox{uc}(\Lambda,n)$.
\begin{df}\label{def:main}
Let $\Lambda\subset \RR^d$ be an integral lattice. 
For $n\in \NN$, the {\sl universally concyclic number }$\mbox{uc}(\Lambda,n)$ (or $\mbox{uc}(n)$ for short) is 
defined by the square of the minimum value among the radii of the hyperspheres that 
pass through exactly $n$ points of $\Lambda$. 
\end{df}
If two lattices $\Lambda_1$ and $\Lambda_2$ are isomorphic, 
then $\mbox{uc}(\Lambda_1,n)=\mbox{uc}(\Lambda_2,n)$ for 
all $n\in \NN$. Therefore, 
$\mbox{uc}(\Lambda,n)$ is an invariant of the lattice $\Lambda$. 
In \cite{M3}, Maehara proposed the following problem: 
\begin{prob}
Determine the $\mbox{uc}(\ZZ^2,n)$ for $n=3,\ldots,10$. 
\end{prob}

In this paper, 
we determine the $\mbox{uc}(\Lambda,n)$ for some $n$ and 
$\Lambda$
whose class number is one. 



The following table provides the computational results.










\newpage 
\begin{thm}\label{thm:main}
Let $K=\QQ(\sqrt{-d})$ as in Theorem \ref{thm:lattice}. Concyclic numbers of two-dimensional Euclidean lattices $\mathcal{O}_K$ of 
class number one for $n\leq 10$ if $d\in\{1,2,3,7,11,19,43,67,163\}$ 
are determined as indicated in Table \ref{table:main}.
\begin{landscape} 
\begin{table}[hbtp]
\caption{}
\label{table:main}
\begin{flushleft}
\begin{tabular}{|c|c|c|c|c|c} 
\noalign{
\hrule 
height0.8pt
}
$-d$  & $d_{K}$ & ${\mathcal{O}_K}$&$\mbox{uc}(3)$&$\mbox{uc}(4)$&$\mbox{uc}(5)$\\ \hline 
\hline
$-1$  & $-2^2$ & $[1, \sqrt{-1}]$ & 
$5^2/2\cdot3^2$ &$1/2$&$5^4/2\cdot3^2$\\ \hline

$-2$  & $-2^3$ & $[1, \sqrt{-2}]$ & 
$3^2/2^3$& $3/2^2$&$3^4/2^3$\\ \hline

$-3$  & $-3$ & $[1, (1+\sqrt{-3})/2]$ & 
$1/3$ &$7/2^2$&$7^2\cdot13^2/11^2$\\ \hline

$-7$  & $-7$ & $[1, (1+\sqrt{-7})/2]$ & 
$2^2/7$ &$2^3/7$&$2^4/7$\\ \hline

$-11$  & $-11$ & $[1, (1+\sqrt{-11})/2]$ & 
$3^2/11$ &$3\cdot5/11$&$3^4/11$\\ \hline

$-19$  & $-19$ & $[1, (1+\sqrt{-19})/2]$ & 
$5^2/19$ &$5\cdot7/19$&$5^4/19$\\ \hline

$-43$  & $-43$ & $[1, (1+\sqrt{-43})/2]$ & 
$11^2/43$ &$11\cdot13/43$&$11\cdot13\cdot17\cdot23/2^2\cdot43$\\ \hline

$-67$  & $-67$ & $[1, (1+\sqrt{-67})/2]$ & 
$17^2/67$ &$17\cdot19/67$&$17\cdot19\cdot23\cdot29/2^2\cdot67$\\ \hline

$-163$ & $-163$ & $[1, (1+\sqrt{-163})/2]$ & 
$41^2/163$ &
$41\cdot43/163$&
$43^2\cdot61^2/3^2\cdot163$\\ 
\noalign{\hrule height0.8pt}
\end{tabular}
\begin{tabular}{c|c|c|c|c|} 
\noalign{\hrule height0.8pt}
\mbox{uc}(6)&$\mbox{uc}(7)$&$\mbox{uc}(8)$&$\mbox{uc}(9)$&$\mbox{uc}(10)$\\ \hline 
\hline
$5^2/2^2$&$5^4\cdot13\cdot17/2\cdot11^2$&$5/2$&$5^2\cdot13^2/2\cdot3^2$
&$5^4/2^2$\\ \hline 

$3^2/2^2$&$3^6/2^3$&$3^3/2^2$&$3^2\cdot11^2/2^3$
&$3^4/2^2$\\ \hline

$1$&$7^2\cdot13\cdot19\cdot43/3\cdot11^2$&$7\cdot13/2^2$&$7^2/3$
&$7^4/2^2$\\ \hline

$2^2$&$2^6/7$&$2^3$&$2^8/7$
&$2^4$\\ \hline

$3^2\cdot5/11$&$3^6/11$&$3^3\cdot5/11$&$3^2\cdot5^2/11$
&$3^4\cdot 5/11$\\ \hline

$5^2\cdot7/19$&$5\cdot7^2\cdot11\cdot17/3^2\cdot19$&$5\cdot7\cdot11/19$&$5^2\cdot7^2/19$
&$5^4\cdot 7/19$\\ \hline

$11^2\cdot13/43$&$11\cdot13^2\cdot17\cdot23/3^2\cdot43$&$11\cdot13\cdot17/43$&$11^2\cdot13^2/43$
&$11^4\cdot 13/43$\\ \hline

$17^2\cdot19/67$&$17\cdot19^2\cdot23\cdot29/3^2\cdot67$&$17\cdot19\cdot23/67$&$17^2\cdot19^2/67$
&$17^4\cdot 19/67$\\ \hline

$41^2\cdot43/163$&                              
$41\cdot43\cdot61\cdot71\cdot83 / 2^2\cdot3^2\cdot163$&
$41\cdot43\cdot47/163$&                        
$41^2\cdot 43^2/163$&                           
$41\cdot47\cdot53\cdot71\cdot83 / 3^2\cdot163$\\

\noalign{\hrule height0.8pt}
\end{tabular}
\end{flushleft}
\end{table}
\end{landscape} 

\end{thm}

We calculated the integer sequences
$\mbox{uc}(\ZZ^2,4n)$ and $\mbox{uc}([1, (1+\sqrt{-3})/2],6n)$
for small $n$,
and speculated that they have simple rules.
Therefore, we have the following problem: 
\begin{prob}\label{conj:main} 
Determine 
$\mbox{uc}(\ZZ^2,4n)$ and 
$\mbox{uc}([1, (1+\sqrt{-3})/2],6n)$ for all $n$. 
\end{prob}
In this paper, we give a partial answer of Problem \ref{conj:main}. 
Namely, we give an exact upperbound of $\mbox{uc}(\ZZ^2, 2^{\ell+2})$ 
and $\mbox{uc}([1, (1+\sqrt{-3})/2],6\cdot 2^m)$. 
\begin{thm}\label{thm:main2} 
Let $\ell$ and $m$ be nonnegative integers, let 
$p_i\ (i=1,2,\ldots)$ be the $i$-th smallest prime 
that is congruent to $1 \pmod{4}$ $($set $p_0:=1$$)$, 
and let $q_j\,(j=1,2,...)$ be the $j$-th smallest prime 
that is congruent to $1 \pmod{3}$ $($set $q_0:=1$$)$. 

\begin{enumerate}

\item [{\rm (1)}]
There exists a circle that passes through exactly $2^{\ell+2}$ points 
$(x,y)$ of $\ZZ^2$: 
\[\left(x-\frac{1}{2}\right)^2+\left(y-\frac{1}{2}\right)^2=\frac{1}{2}\, 
\prod_{k=0}^\ell p_k.\] 
Therefore, we have 
\[
{\rm uc}(\ZZ^2,2^{\ell+2})\leq \frac{1}{2}\prod_{k=0}^\ell p_k. 
\]

\item [{\rm (2)}]
The number of the integer solutions of the following equation
\[
\left(x + y \frac{1 + \sqrt{-3}}{2}\right) 
\left(x + y \frac{1 - \sqrt{-3}}{2}\right)
=x^2+xy+y^2=\prod_{k=0}^m q_k\]
is $6\cdot 2^{m}$. 
This means that 
the circle 
\[
|z|=\left(\prod_{k=0}^m q_k\right)^{\frac{1}{2}}
\]
passes through exactly $6\cdot 2^{m}$ points of 
$[1, (1+\sqrt{-3})/2]$. 
Therefore, we have 
\[
{\rm uc}([1, (1+\sqrt{-3})/2],6\cdot 2^{m})\leq \prod_{k=0}^m q_k. 
\]

\end{enumerate}
\end{thm}

In Section \ref{sec:alg}, 
we give the computational algorithm used in Theorem \ref{thm:main}. 
In Section \ref{sec:proof}, 
we provide the proof of Theorem \ref{thm:main2}. 
In Section \ref{sec:problem}, 
we present further problems. 

All the computer calculations in this paper 
were done by Mathematica \cite{Ma} and C Programming Language \cite{C}.

\section{Algorithm}\label{sec:alg}

In this section, we give the algorithm used to find the square of 
the minimum value among the radii of the hyperspheres 
that pass through exactly $n$ points of $\Lambda$.

Assume that $\Lambda$ is one of $\mathcal{O}_K$ in Theorem \ref{thm:main}. 
Let $\ell$ be a positive integer, and let $R\subset{\Lambda}$ 
be the set of $(x,y)$ that satisfies $x^2+y^2\leq\ell^2$, $y\geq0$ and $y\geq-\sqrt{d}x$ if $d=2,3,7,11,19,43,67,163$
(if $d=1$, then let $R$ be the set of $(x,y)$ that satisfies 
$x^2+y^2\leq\ell^2$, $x\geq0$ and $y\geq0$).
We shall try to create a hypersphere by taking three vertices 
on $R$.
Notice that a hypersphere is determined uniquely 
by taking three vertices over $\Lambda$. 

First, we shall explain how to plot the three vertices 
on $R$. Let $(x_i,y_i)$ be the $i$-th vertex $(i=1,2,3)$. 
Set $(x_1,y_1)=(0,0)$, and let $(x_2,y_2)$ vary such that it plots every vertex 
$(x,y)\in\Lambda$ such that $y/x<\sqrt{d}$ on $R$. Then, we let 
$(x_3,y_3)$ vary such that it plots every vertex $(x,y)\in\Lambda$,
except for $(x,y)\in\Lambda$ such that $y=0$ on $R$. 
This algorithm will provide every hypersphere passing through $(0,0)$ that can be generated 
by any $(x,y)\in\Lambda$ on $R$.

Next, we shall explain how to obtain the coordinates for the center and 
the square of the radius of a hypersphere. 
Let $(x_0,y_0)$ be the center of a hypersphere, and let
$D$ be the square of the radius of the hypersphere.
Then,
\begin{align*}
x_0&=-(y_2y_3^2 + (-y_2^2 - x_2^2)y_3 + x_3^2y_2)/(2x_2y_3 - 2x_3y_2),\\
y_0&=(x_2y_3^2-x_3y_2^2+x_2x_3^2-x_2^2x_3)/(2x_2y_3-2x_3y_2),\\
D&=\left(\frac{\sqrt{y_2^2+x_2^2}\sqrt{\alpha}}{2x_2y_3-2x_3y_2}\right)^2,
\end{align*}
where $\alpha=y_3^4-2y_2y_3^3+y_2^2y_3^2+2x_3^2y_3^2-2x_2x_3y_3^2+x_2^2y_3^2-2x_3^2y_2y_3+x_3^2y_2^2+x_3^4-2x_2x_3^3+x_2^2x_3^2$. 

Next, we explain how to enumerate the number of lattice points 
$(x,y)\in\Lambda$ such that $(x-x_0)^2+(y-y_0)^2=D$. 
Let $x_4\in\Lambda$ move from $[x_0-\sqrt{D}]$ to $[x_0+\sqrt{D}]+1$ 
, where $[\ ]$ is the Gauss symbol.
For the equation $(x_4-x_0)^2+(y_4-y_0)^2=D$, solve for $y_4$: 
$y_4=y_0\pm\sqrt{-x_4^2+2x_0x_4+D-x_0^2}$.
Set $c_p=0$.
If $x_4\equiv0 \pmod{1}$ and $y_4\equiv0 \pmod{\sqrt{d}}$, or
if $x_4\equiv 1/2 \pmod{1}$ and $y_4\equiv\sqrt{d}/2 \pmod{\sqrt{d}}$,
then $c_p=c_p+1$ (in the case of $d=3,7,11,19,43,67,163$).
If $x_4\equiv0 \pmod{1}$ and $y_4\equiv0 \pmod{\sqrt{d}}$,
then $c_p=c_p+1$ (in the case of $d=1,2$).
It is seen that $c_p$ denotes the number of 
lattice points $(x,y)\in\Lambda$ such that 
$(x-x_0)^2+(y-y_0)^2=D$ after moving $x_4$ 
from $[x_0-\sqrt{D}]$ to $[x_0+\sqrt{D}]+1$. 
Therefore, we can obtain the hypersphere 
that passes through exactly $c_p$ points.
 
Using the above method, since we can find the hyperspheres
that pass through exactly $c_p$ points for any $n\in\NN$, 
we can obtain 
the square of the minimum value of the radius by selecting the smallest radius of any of the hyperspheres that pass through exactly 
$n$ points of $\Lambda$.

\section{Proof of Theorem \ref{thm:main2}}\label{sec:proof}


First, 
we claim that the circle $(2x-1)^2+(2y-1)^2=2\prod_{k=0}^\ell p_k$ 
passes through exactly $2^{\ell+2}$ points of $\ZZ^2$.
By Fermat's $4n+1$ Theorem, 
for all $p_j$, there exists $a_j$,$b_j\in\ZZ$ such that 
$p_j=a_j^2+b_j^2$.
Therefore, $p_j=(a_j+ib_j)(a_j-ib_j)$.
Notice that $a_j+ib_j$ and $a_j-ib_j$ are irreducible elements 
over $\ZZ[i]$. 
Since $2=1^2+1^2$, $2=(1+i)(1-i)$.
Hence $\omega\overline{\omega}=2\,\prod_{k=0}^\ell p_k
=(1+i)(1-i)\prod_{k=0}^\ell p_k
=(1+i)(1-i)\prod_{k=0}^\ell (a_k+ib_k)(a_k-ib_k)$,
where $\omega\in\ZZ[i]$.
We consider the number of possible outcomes for  $\omega$. We can express $\omega$ as follows:\\
$\omega=u(1+i)^{\epsilon_0}(1-i)^{1-\epsilon_0}(a_1+ib_1)^{\epsilon_1}(a_1-ib_1)^{1-\epsilon_1}\cdots (a_\ell+ib_\ell)^{\epsilon_\ell}(a_\ell-ib_\ell)^{1-\epsilon_\ell}$, where 
$u=\pm1,\pm{i}$, and $\epsilon_n=0,1$ $(n=0,1,\ldots,\ell)$. 

It is easily seen that the choice of $(1+i)$ or $(1-i)$ does not 
depend on the number of possible outcomes of 
$\omega$, since the absolute value of the real part and 
the imaginary part of $(1+i)$ and $(1-i)$ is the same.

Consequently, the number of possible outcomes of $\omega$ is 
$4\cdot 2^{\ell+1}/2=2^{\ell+2}$ over $\ZZ[i]$.
From this, the number of $(X,Y)\in\ZZ^2$ such that $X^2+Y^2=2\,\prod_{k=0}^\ell p_k$ is $2^{\ell+2}$.

Next, we claim that they all correspond to the lattice point $(x,y)\in\ZZ^2$ such that $(2x-1)^2+(2y-1)^2=2\,\prod_{k=0}^\ell p_k$.
Since $X^2,Y^2 \equiv 0,1 \pmod{4}$ and
$2\prod_{k=0}^\ell p_k \equiv 2 \pmod{4}$,
$X^2+Y^2=\prod_{k=0}^\ell p_k$ implies that $X^2\equiv1$ and $Y^2\equiv1$$\pmod{4}$. Moreover, it implies that $X \equiv 1$ and $Y \equiv 1$ $\pmod{2}$.
Therefore, the number of lattice points $(x,y)\in\ZZ^2$
such that $(2x-1)^2+(2y-1)^2=2\,\prod_{k=0}^l p_k$ is equivalent to
the number of $(X,Y)\in\ZZ^2$ such that  $X^2+Y^2=2\,\prod_{k=0}^l p_k$, 
$2x-1\equiv-1\equiv1(\equiv X)$ $\pmod{2}$ and 
$2y-1\equiv-1\equiv1(\equiv Y)$ $\pmod{2}$.

Thus, the number of lattice points $(x,y)\in\ZZ^2$ such that $(2x-1)^2+(2y-1)^2=2\,\prod_{k=0}^\ell p_k$ is just $2^{\ell+2}$.


Next, we claim that the number of the integer solutions of 
the following equation
\[
\left(x + y \frac{1 + \sqrt{-3}}{2}\right) 
\left(x + y \frac{1 - \sqrt{-3}}{2}\right)
=x^2+xy+y^2=\prod_{k=0}^m q_k\]
is $6\cdot 2^{m}$. 

The proof is similar to the first part. 
Set $\zeta=1+\sqrt{-3}/2$.
Then, for all $q_i$, there exists $a_i$,$b_i\in\ZZ$ 
such that 
\begin{align*}q_i&=a_i^2+a_ib_i+b_i^2=
\left(a_i+\frac{b_i}{2}+\frac{\sqrt{-3}b_i}{2}\right)
\left(a_i+\frac{b_i}{2}-\frac{\sqrt{-3}b_i}{2}\right)\\
&=(a_i+b_i\zeta)(a_i+b_i\overline\zeta).
\end{align*}
Notice that $a_i+b_i\zeta$ and $a_i+b_i\overline\zeta$ 
are irreducible elements over $\ZZ[\zeta]$, and
$\tau\overline{\tau}=\prod_{k=0}^m q_k
=\prod_{k=0}^m (a_k+b_k\zeta)(a_k+b_k\overline\zeta)$, where $\tau\in\ZZ[\zeta]$.
We consider the number of possible outcomes for $\tau$.
We can express $\tau$ as follows:
$\tau=u(a_1+b_1\zeta)^{\mu_1}(a_1+b_1\overline{\zeta})^{1-\mu_1}
\cdots(a_m+b_m\zeta)^{\mu_m}(a_m+b_m\overline{\zeta})^{1-\mu_m}$, 
where $u=\pm1,\pm\zeta,\pm\overline\zeta$, $\mu_n=0,1$ 
$(n=1,\ldots,m)$. 

As a consequence, 
the number of possible outcomes of $\tau$ is $6\cdot 2^m$ over 
$\ZZ[\zeta]$.
From this, the number of $x+y\zeta\in\ZZ[\zeta]$ such that 
$x^2+xy+y^2=\,\prod_{k=0}^m q_k$ is $6\cdot 2^m$. Now, 
since it can be seen that $\ZZ[\zeta]$ is equivalent to 
$[1, (1+\sqrt{-3})/2]$, 
the circle 
\[
|z|=\left(\prod_{k=0}^m q_k\right)^{\frac{1}{2}}
\]
passes through exactly $6\cdot 2^{m}$ points of 
$[1, (1+\sqrt{-3})/2]$. 

\begin{rem}
We remark that 
the conditions in Theorem \ref{thm:main2} 
``the $i$-th smallest prime" and ``the $j$-th smallest prime"
do not use in the proof of Theorem \ref{thm:main2}. 
For example, the number of solutions (points of $\ZZ^2$) 
is determined by the number of primes appearing in the product
\[
\prod _{k=0}^\ell p_k. 
\]
On the other hand, we need these conditions in order to answer Problem \ref{conj:main}. 
\end{rem}

\section{Further problems}\label{sec:problem}

\begin{enumerate}

\item [(1)]
Find a law in the table of Theorem \ref{thm:main}.

\item [(2)]
For $n=\{3,\ldots,10\}$, determine the $\mbox{uc}(\Lambda,n)$ for 
$\Lambda=\ZZ^3$ and $\ZZ^4$. 

\item [(3)]
Let
\begin{align*}
\left\{
\begin{array}{l}
{\bf e_1} = \displaystyle \frac{1}{\sqrt{12}} (1, 0, 0, 0) \\
{\bf e_2} = \displaystyle \frac{7}{\sqrt{12}} (0, 1, 0, 0) \\
{\bf e_3} = \displaystyle \frac{13}{\sqrt{12}} (0, 0, 1, 0) \\
{\bf e_4} = \displaystyle \frac{19}{\sqrt{12}} (0, 0, 0, 1). 
\end{array}
\right.
\end{align*}
Then, we define the two lattices, 
$L_1:=\langle {\bf u_i}\mid i=1,\ldots,4\rangle$ and 
$L_2:=\langle {\bf v_i}\mid i=1,\ldots,4\rangle$, 
where 
\begin{align*}
\left\{
\begin{array}{l}
{\bf u_1} = 3 {\bf e_1} - {\bf e_2} - {\bf e_3} - {\bf e_4} \\
{\bf u_2} = {\bf e_1} + 3 {\bf e_2} + {\bf e_3} - {\bf e_4} \\
{\bf u_3} = {\bf e_1} - {\bf e_2} + 3 {\bf e_3} + {\bf e_4} \\
{\bf u_4} = {\bf e_1} + {\bf e_2} - {\bf e_3} + 3 {\bf e_4}, \\ 
\end{array}
\right.
\left\{
\begin{array}{l}
{\bf v_1} = -3 {\bf e_1} - {\bf e_2} - {\bf e_3} - {\bf e_4} \\
{\bf v_2} = {\bf e_1} - 3 {\bf e_2} + {\bf e_3} - {\bf e_4} \\
{\bf v_3} = {\bf e_1} - {\bf e_2} - 3 {\bf e_3} + {\bf e_4} \\
{\bf v_4} = {\bf e_1} + {\bf e_2} - {\bf e_3} - 3 {\bf e_4}. 
\end{array}
\right.
\end{align*}

In \cite{CS}, it was shown that 
the theta series of $L_1$ and $L_2$ are the same, 
namely, 
the number of lattice vectors of norm $m$ are the same for all $m$.
However, these two lattices are nonisomorphic, and 
the proof of this fact is not easy \cite{CS}. 

Therefore, we have the following problem: 
Determine the $\mbox{uc}(L_1,n)$ and $\mbox{uc}(L_2,n)$ 
for some $n$, and show that $L_1$ and $L_2$ are nonisomorphic.

\end{enumerate}
\section*{Acknowledgments}
The authors would like to thank the anonymous
reviewers for their beneficial comments 
on an earlier version of the manuscript. 
This work was supported by JSPS KAKENHI (18K03217).

\end{document}